\documentclass[12pt]{amsart}

\usepackage{amssymb,amsmath,cite,graphicx}
\usepackage[sc]{mathpazo}

\newcommand\BS{\boldsymbol}
\newcommand\dif{\mathrm{d}}
\newcommand\deriv[2]{\frac{\dif #1}{\dif #2}}
\newcommand\parderiv[2]{\frac{\partial #1}{\partial #2}}

\numberwithin{equation}{section}

\begin{document}

\title[Linear response of the Lyapunov exponent]{Linear response of
  the Lyapunov exponent to a small constant perturbation}

\address{Department of Mathematics, Statistics and Computer
  Science\\University of Illinois at Chicago\\851 S. Morgan st. (M/C
  249)\\ Chicago, IL 60607}

\author{Rafail V. Abramov}

\email{abramov@math.uic.edu}


\date{\today}

\pagestyle{myheadings}

\begin{abstract}
In the current work we demonstrate the principal possibility of
prediction of the response of the largest Lyapunov exponent of a
chaotic dynamical system to a small constant forcing perturbation via
a linearized relation, which is computed entirely from the unperturbed
dynamics. We derive the formal representation of the corresponding
linear response operator, which involves the (computationally
infeasible) infinite time limit. We then compute suitable finite-time
approximations of the corresponding linear response operator, and
compare its response predictions with actual, directly perturbed and
measured, responses of the largest Lyapunov exponent. The test
dynamical system is a 20-variable Lorenz 96 model, run in weakly,
moderately, and strongly chaotic regimes. We observe that the
linearized response prediction is a good approximation for the
moderately and strongly chaotic regimes, and less so in the weakly
chaotic regime due to intrinsic nonlinearity in the response of the
Lyapunov exponent, which the linearized approximation is incapable of
following.
\end{abstract}

\maketitle

\section{Introduction}
\label{sec:intro}

The largest Lyapunov exponent (henceforth {\em the} Lyapunov exponent)
is the cornerstone measure of chaos and uncertainty in complex
nonlinear dynamics \cite{Ose,EckRue,Ott,Str}.  It shows the average
exponential rate of divergence (if positive) or convergence (if
negative) of two nearby solutions of a dynamical system, usually a
system of nonlinear differential equations of first order, sometimes
with stochastic forcing. Practical methods of computing the Lyapunov
exponents from the long-term time series of a dynamical system have
been developed in \cite{WolSwiSwiVas,Kan}.

In this work we develop a linear approximation to the response of the
Lyapunov exponent to a small constant perturbation of the
corresponding dynamical system. This linear approximation is computed
from the long-term time series of the unperturbed dynamics, and, in a
way, is a generalization of the well-known Fluctuation-Dissipation
theorem \cite{MajAbrGro,Ris,Web}. The advantage of the approach is
that it effectively computes all possible responses to all possible
perturbations in the corresponding phase space at once, as well as
provides a convenient framework for the inverse problem, where the
perturbation has to be computed to satisfy the given response. The
approach we develop here is not generally restricted to constant
perturbations, and can be adapted to arbitrary differentiable
perturbations of vector fields, provided that the Lyapunov exponent
varies sufficiently smoothly under the perturbations.

We consider a system of autonomous nonlinear ordinary differential
equations of the form
\begin{equation}
\label{eq:dyn_sys}
\deriv{\BS x}t = \BS f(\BS x),
\end{equation}
where $\BS x=\BS x(t)$ is a $\mathbb R^N$-vector valued bounded
function of time $t$ for some positive integer $N$, representing the
unknown time-dependent solution of the system, and $\BS f(\BS x)$ is a
nonlinear differentiable vector field, $\BS f:\mathbb R^N\to\mathbb
R^N$.

The largest Lyapunov exponent describes the average exponential rate
of separation of solutions $\BS x(t)$ and $\BS y(t)$ of
\eqref{eq:dyn_sys} for two nearby initial conditions $\BS x_0$ and $\BS
y_0$, in the infinite time limit:
\begin{equation}
\label{eq:lyap_exp}
\lambda=\lim_{t\to\infty}\lim_{\BS y_0\to\BS x_0}\frac
1t\ln\frac{\|\BS y(t)-\BS x(t)\|}{\|\BS y_0-\BS x_0\|}.
\end{equation}
Above, $\|\BS x\|$ denotes the usual Euclidean norm of $\BS x$. The
famous Oseledec's multiplicative ergodic theorem \cite{Ose} states
that for almost all, in the appropriate probabilistic sense, starting
conditions $\BS x_0$, the limit in \eqref{eq:lyap_exp} converges to
the same value independently of choice of the initial condition $\BS
x_0$. The relation in \eqref{eq:lyap_exp} means that, on average, the
distance between two nearby solutions $\BS x(t)$ and $\BS y(t)$ can be
estimated as
\begin{equation}
\|\BS y(t)-\BS x(t)\|~\sim e^{\lambda t}\|\BS y_0-\BS x_0\|.
\end{equation}
Thus, if $\lambda$ is positive, almost any two nearby solutions of
\eqref{eq:dyn_sys} diverge from each other exponentially rapidly in
time, with $\lambda$ being the average exponential rate of
divergence. In this case, the dynamical system in \eqref{eq:dyn_sys}
is said to be {\em chaotic} \cite{Ott,Str}. Chaotic dynamical systems
can be encountered in fluid dynamics \cite{Are}, turbulence
\cite{Kra2,Rue0} and geophysical science \cite{Lor2,Kal}, and are the
subject of our current study.

There is more simple formula for the Lyapunov exponent which excludes
the spatial limit.  Observe that, for two nearby solutions $\BS x(t)$
and $\BS y(t)$, the difference between them can be approximated as
\begin{equation}
\begin{array}{c}
\displaystyle\deriv{\BS v}t = D\BS f(\BS x)\BS v+o(\|\BS
v\|),\\ \displaystyle\BS v(t)=\BS y(t)-\BS x(t),
\end{array}
\end{equation}
where $D\BS f$ is the Jacobian of $\BS f$ (the matrix of partial
derivatives of $\BS f$). As $\|\BS v\|\to 0$, the term $o(\|\BS
v\|)$ above becomes negligible in comparison with the rest of the
terms, which results in the linear equation for $\BS v(t)$ in the
limit, with $\BS x(t)$ computed in parallel from
\eqref{eq:dyn_sys}:
\begin{equation}
\label{eq:dyn_sys_tangent}
\deriv{\BS v}t = D\BS f(\BS x)\BS v,\qquad
\lambda=\lim_{t\to\infty}\frac 1t\ln\frac{\|\BS v(t)\|}{\|\BS v_0\|}.
\end{equation}
Here, $\BS v(t)$ is a {\em tangent vector}, with $\|\BS v_0\|=1$, by
convention. Here, $\BS v(t)$ does not have to be small, since it can
be scaled by an arbitrary constant factor due to the linearity and
cancellation in \eqref{eq:dyn_sys_tangent}. For chaotic systems,
$\|\BS v(t)\|$ grows exponentially fast, so we avoid the resulting
numerical instability with periodic renormalization of $\BS v(t)$ and
the corresponding adjustment to $\lambda$ for compensation (for
details, see, for example, \cite{EckRue}).

One can also replace the infinite time limit in
\eqref{eq:dyn_sys_tangent} with a time average (or ensemble average,
via Birkhoff's theorem \cite{Bir}). First, factor the tangent vector
$\BS v$ from \eqref{eq:dyn_sys_tangent} into the product of its norm
$\|\BS v\|$, and unit vector $\BS w=\BS v/\|\BS v\|$:
\begin{equation}
\label{eq:polar}
\BS v=\|\BS v\|\BS w.
\end{equation}
Then, for the unperturbed system in \eqref{eq:dyn_sys}, the equation
for $\|\BS v(t)\|$ is easily derived as
\begin{equation}
\label{eq:norm_equation}
\deriv{\ln\|\BS v\|}t=\BS w^TD\BS f(\BS x)\BS w.
\end{equation}
The substitution of \eqref{eq:polar} and \eqref{eq:norm_equation} into
the formula for the Lyapunov exponent in \eqref{eq:dyn_sys_tangent}
yields, with $n_0=1$,
\begin{equation}
\label{eq:lyap_exp_time_average}
\begin{split}
\lambda=\lim_{t\to\infty}\frac 1t\ln\|\BS
v(t)\|=\lim_{t\to\infty}\frac 1t\int_0^t\deriv{\ln\|\BS v(s)\|}s\dif s=\\=
\lim_{t\to\infty}\frac 1t\int_0^t\BS w^T(s)D\BS f(\BS x(s))\BS
w(s)\dif s.
\end{split}
\end{equation}
Optionally, one can invoke Birkhoff's ergodic theorem \cite{Bir} and
obtain the statistical average
\begin{equation}
\label{eq:lyap_exp_average}
\lambda=\int\BS w^TD\BS f(\BS x)\BS w\,\dif\rho(\BS x,\BS w),
\end{equation}
where $\rho(\BS x,\BS w)$ is the joint invariant distribution measure
of $\BS x$ and $\BS w$. The average formulas in
\eqref{eq:lyap_exp_time_average} and \eqref{eq:lyap_exp_average} is
what we need to derive the linear response approximation formula for
the Lyapunov exponent below.

The manuscript is organized as follows: in Section \ref{sec:response}
we present the formal derivation of the linear response operator for
the Lyapunov exponent, in Section \ref{sec:computation} we show the
results of numerical computational tests of the derived formula with a
chaotic nonlinear test model, and in Section \ref{sec:conclusions} we
summarize the results of the work. Additionally, Appendix
\ref{sec:app1} contains technical details of derivations for Section
\ref{sec:response}, while Appendix \ref{sec:app2} outlines the
computational discretization of the response formula.

\section{Linear response to small constant perturbation}
\label{sec:response}

One can introduce different types of small perturbations in the
right-hand side of \eqref{eq:dyn_sys}, and the approach, developed
below, can be adapted to an arbitrary differentiable perturbation
vector field. However, for simplicity of presentation, here we
consider the most basic case of a constant vector perturbation $\BS
p$:
\begin{equation}
\label{eq:dyn_sys_pert}
\deriv{\BS x}t = \BS f(\BS x)+\BS p.
\end{equation}
While this perturbation may seem trivial, there is no easy way to tell
in general how the dynamical properties of \eqref{eq:dyn_sys} respond
to such a perturbation (e.g. complex bifurcations of the resulting
flow may occur, fixed points/periodic orbits created/destroyed, etc).
Even if no bifurcations occur (which is what we assume here throughout
the work), the solution $\BS x(t)$ generally changes nonlinearly under
finite constant perturbations, and so should do the Lyapunov exponent
in \eqref{eq:lyap_exp_time_average} and \eqref{eq:lyap_exp_average}.

Below we develop a linear approximation of the response of the
Lyapunov exponent to the constant perturbation $\BS p$ in
\eqref{eq:dyn_sys_pert}, under the condition that $\|\BS p\|$ is
sufficiently small. For the perturbed system in
\eqref{eq:dyn_sys_pert}, the corresponding formula for the Lyapunov
exponent reads
\begin{equation}
\label{eq:lambda_pert_1}
\lambda_{\BS p}=\int\BS w^TD\BS f(\BS x)\BS w\,\dif\rho_{\BS p}(\BS
x,\BS w),
\end{equation}
where $\rho_{\BS p}(\BS x,\BS w)$ is the joint invariant distribution
of $\BS x$ and $\BS w$ for the perturbed system in
\eqref{eq:dyn_sys_pert}, with a shorthand notation $\rho_{\BS 0}(\BS
x,\BS w)=\rho(\BS x,\BS w)$. Here we adopt the flow notations $\BS
x(t)=\phi_{\BS p}^t\BS x$, $\BS w(t)=\psi_{\BS p,\BS x}^t\BS w$, where
subscripts are dropped from the initial conditions $\BS x_0$ and $\BS
w_0$, and denote $\phi_{\BS 0}^t=\phi^t$, $\psi_{\BS 0,\BS
  x}^t=\psi_{\BS x}^t$. The next step is to represent $\rho_{\BS
  p}(\BS x,\BS w)$ as the pushforward measure \cite{Rue2} of $\rho(\BS
x,\BS w)$:
\begin{equation}
\begin{split}
\rho_{\BS p}(\BS x,\BS w)=\lim_{t\to\infty}\rho(\phi_{\BS p}^{-t}\BS
x,\psi_{\BS p,\BS x}^{-t}\BS w),
\end{split}
\end{equation}
which yields, for $\lambda_{\BS p}$,
\begin{equation}
\begin{split}
\lambda_{\BS p}=\lim_{t\to\infty}\int\BS w^TD\BS f(\BS x)\BS
w\,\dif\rho(\phi_{\BS p}^{-t}\BS x,\psi_{\BS p,\BS x}^{-t}\BS
w)=\\=\lim_{t\to\infty}\int\left(\psi_{\BS p,\BS x}^t\BS w\right)^TD\BS
f(\phi_{\BS p}^t\BS x)\psi_{\BS p,\BS x}^t\BS w\,\dif\rho(\BS x,\BS w),
\end{split}
\end{equation}
where the last equality follows from the change of variables $\BS
x\to\phi_{\BS p}^t\BS x$, $\BS w\to\psi_{\BS p,\BS x}^t\BS w$. Next,
we use the invariance of the measure $\rho(\BS x,\BS w)$ with respect
to the unperturbed flows $\phi^t\BS x$ and $\psi_{\BS x}^t\BS w$ (that
is, $\rho(\BS x,\BS w)=\rho(\phi^{-t}\BS x,\psi_{\BS x}^{-t}\BS w)$
for any $t$) and rewrite the above relation equivalently as
\begin{subequations}
\label{eq:lyap_exp_average_pert}
\begin{equation}
\lambda_{\BS p}=\lim_{t\to\infty}\int L_{\BS p}\,\dif\rho(\BS x,\BS
w),
\end{equation}
\begin{equation}
\begin{split}
L_{\BS p}=\big(\psi_{\BS p,\phi^{-t}\BS x}^t\psi_{\BS x}^{-t}\BS
w\big)^TD\BS f(\phi_{\BS p}^t\phi^{-t}\BS x)\psi_{\BS p, \phi^{-t}\BS
  x}^t\psi_{\BS x}^{-t}\BS w.
\end{split}
\end{equation}
\end{subequations}
Our next step is to approximate the difference $\lambda_{\BS p}-\lambda$ via
the linear relation
\begin{equation}
\label{eq:lin_rel}
\lambda_{\BS p}-\lambda \approx
\lim_{t\to\infty}\left(\int\left.\parderiv{L_{\BS p}}{\BS
  p}\right|_{\BS p=\BS 0}\,\dif\rho(\BS x,\BS w)\right)\cdot\BS p,
\end{equation}
under the assumption that the derivative exists in the infinite time
limit (so-called {\em structural stability} \cite{AlbSpr,Man}).
Usually, it is the case with deterministic dynamical systems with
Sinai-Ruelle-Bowen invariant measures \cite{Sin,Rue3,You}, and many
random dynamical systems, including It\^o diffusions \cite{Arn,Kun}.
After some calculations (which are provided in Appendix
\ref{sec:app1}) we obtain the linear relation in the form of the
fluctuation-response time correlation functions of the unperturbed
system in \eqref{eq:dyn_sys}:
\begin{subequations}
\label{eq:corr}
\begin{equation}
\lambda_{\BS p}-\lambda\approx\lim_{t\to\infty}\BS r(t)\cdot\BS p,
\end{equation}
\begin{equation}
\BS r(t)=\int_0^t\BS c_1(\tau)\dif\tau+\int_0^t\dif\tau\int_\tau^t\BS
c_2(\tau,s)\dif s,
\end{equation}
\begin{equation}
\begin{split}
\BS c_1(\tau)=\lim_{T\to\infty}\frac 1T\int_0^T\BS w(t^\prime)^TD^2\BS
f(\BS x(t^\prime)):\\:\big(\BS w(t^\prime)\otimes\BS T_{\BS
  x(t^\prime-\tau)}^\tau\big)\dif t^\prime,
\end{split}
\end{equation}
\begin{equation}
\begin{split}
\BS c_2(\tau,s)=\lim_{T\to\infty}\frac 1T\int_0^T\BS
w^T(t^\prime)\big(D\BS f(\BS x(t^\prime))^T+\\+D\BS f(\BS
x(t^\prime))\big)\big(\BS I-\BS w(t^\prime)\BS w^T(t^\prime)\big)\BS
T_{\BS x(t^\prime-\tau)}^\tau\times\\\times D^2\BS f(\BS
x(t^\prime-\tau)):\Big(\BS T_{\BS x(t^\prime)}^{-\tau}\BS
w(t^\prime)\otimes\BS T_{\BS x(t^\prime-s)}^{s-\tau}\Big)\dif
t^\prime.
\end{split}
\end{equation}
\end{subequations}
Above, the Frobenius product ``$:$'' is computed over the two
covariant indices of $D^2\BS f$. The $N\times N$ matrix $\BS T_{\BS
  x}^t$ is the {\em tangent map} of $\phi^t\BS x$:
\begin{equation}
\BS T_{\BS x}^t=\parderiv{}{\BS x}\phi^t\BS x.
\end{equation}
The equation, which is used for the computation of the tangent map, is
given in \eqref{eq:tangent_map}.

\section{Computational approximation and numerical testing}
\label{sec:computation}

From the formula in \eqref{eq:corr} it follows that the linear
approximation of the response of the Lyapunov exponent is the integral
over a specially crafted time-lag correlation function with the
infinite upper limit of integration. In general, the computation of
such limit is feasible, as long as the integrand (the time-lag
correlation under the integral) decays sufficiently rapidly to zero
with increasing time, so that the integral over it could be truncated
to some finite upper limit. However, for a chaotic and mixing
dynamical system in \eqref{eq:dyn_sys}, the situation is complicated
by the fact that the tangent map $\BS T_{\BS x}^t$ of the unperturbed
system in \eqref{eq:dyn_sys} grows exponentially rapidly in $t$ for
almost any $\BS x$, with the exponential rate of $\lambda$, even
though its ensemble average decays in $t$ \cite{EckRue,Rue2}. This
causes the numerical instability in the form of precision loss; many
extermely large numbers must add up to small numbers, which does not
happen in the finite precision computer arithmetic
\cite{Abr5,Abr6,Abr7,AbrMaj4,AbrMaj5,AbrMaj6}.

\begin{figure}
\includegraphics[width=0.5\textwidth]{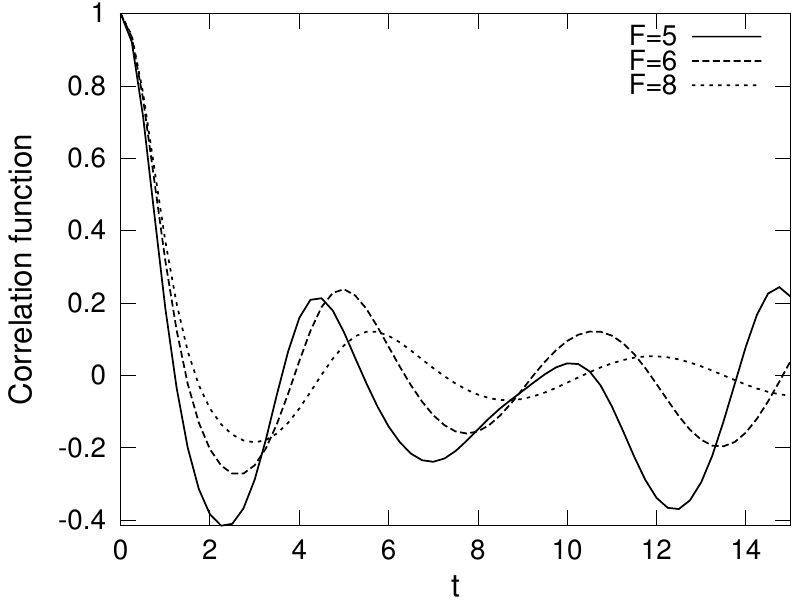}
\caption{Time-lag autocorrelation functions of the 20-variable Lorenz
  96 model for the regimes $F=5,6,8$.}
\label{fig:l96_corr}
\end{figure}

\begin{figure}
\includegraphics[width=\textwidth]{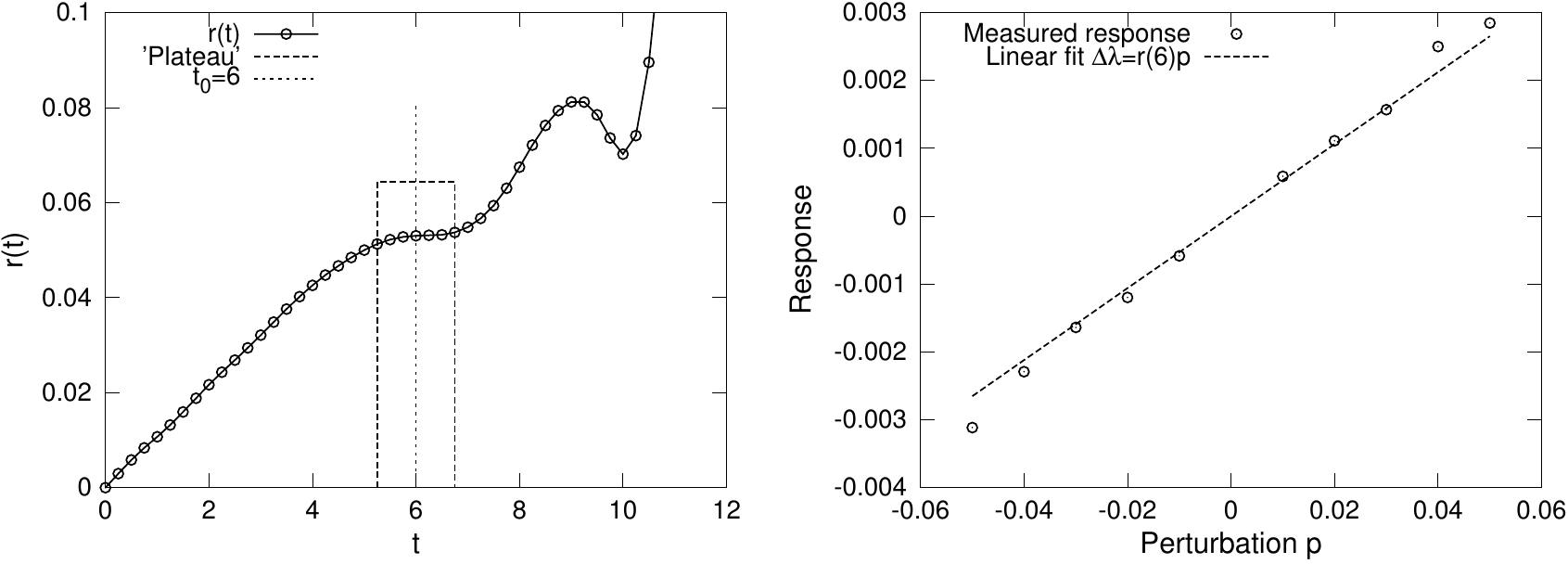}
\caption{Linear response of the Lyapunov exponent for the 20-variable
  Lorenz 96 model \eqref{eq:l96} with $F=8$ ($\lambda=0.4253$). Left:
  the plot $\BS r(t)$ as a function of finite response time, and the
  ``plateau'' of saturated response before manifestation of numerical
  instability. The time $t_0=6$ is chosen arbitrarily from within the
  ``plateau'' window. Right: the measured response of the Lyapunov
  exponent (via direct perturbations) vs the linear fit $\BS
  r(t_0)\cdot\BS p$.}
\label{fig:l96_F8}
\end{figure}

As a result, we are forced to consider a finite-time response
approximation
\begin{equation}
\label{eq:finite_resp}
\Delta\lambda\approx\BS r(t_0)\cdot\BS p,
\end{equation}
for a finite response time $t_0$, rather than an infinite time limit.
The range of $t$, for which the response operator $\BS r(t)$ is
practically computable without significant numerical instability, is
usually proportional to the {\em e-folding time} $\lambda^{-1}$
\cite{AbrMaj4,AbrMaj5,AbrMaj6}. The finite response time $t_0$ for
\eqref{eq:finite_resp} will have to be chosen from that range.
\begin{figure}
\includegraphics[width=\textwidth]{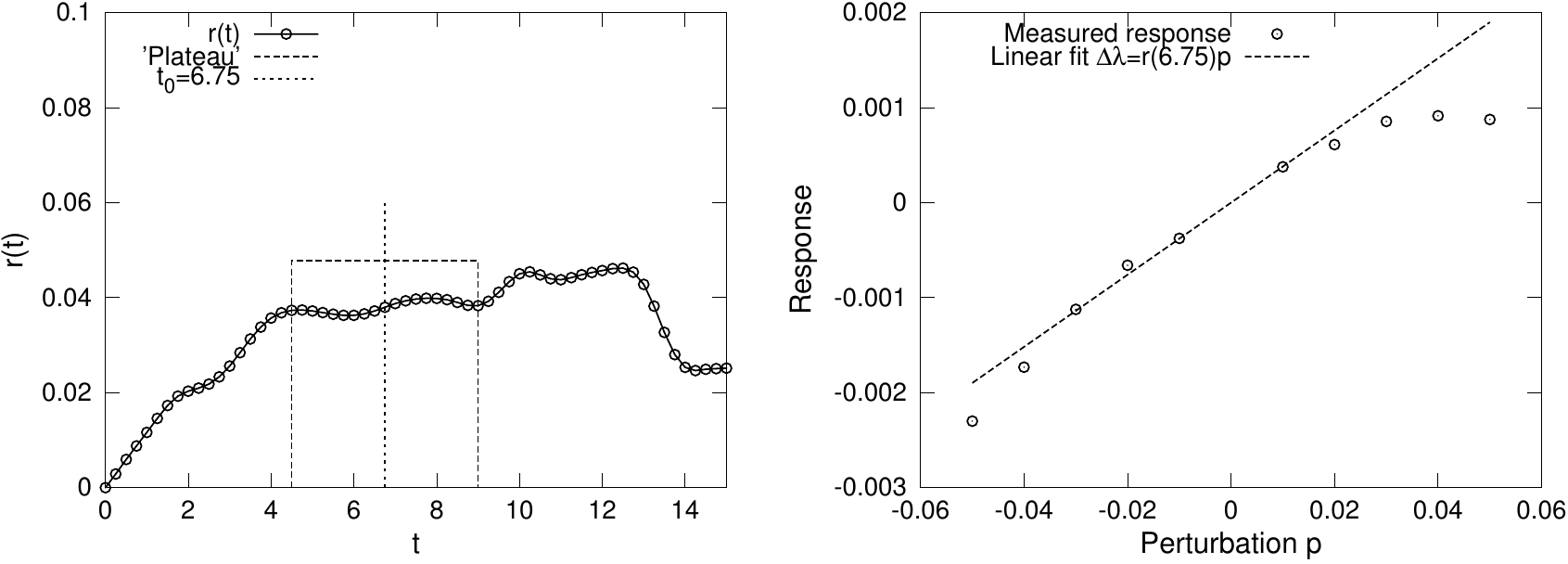}
\caption{Linear response of the Lyapunov exponent for the 20-variable
  Lorenz 96 model \eqref{eq:l96} with $F=5$ ($\lambda=0.2265$). Left:
  the plot $\BS r(t)$ as a function of finite response time, and the
  ``plateau'' of saturated response before manifestation of numerical
  instability. The time $t_0=6.75$ is chosen arbitrarily from within
  the ``plateau'' window. Right: the measured response of the Lyapunov
  exponent (via direct perturbations) vs the linear fit $\BS
  r(t_0)\cdot\BS p$.}
\label{fig:l96_F5}
\end{figure}
As the test system, we consider the rescaled Lorenz 96 model
\cite{Abr5,Abr6,Abr7,Abr8,Abr9,Abr10,Abr11,AbrMaj4,AbrMaj5,AbrMaj6,MajAbrGro,Lor,LorEma},
which is a simple nonlinear chaotic forced-dissipative system with a
band of linearly unstable waves with oppositely directed phase and
group velocities, similar to the Rossby waves in the midlatitudinal
troposphere. The rescaled Lorenz 96 model is given by the system of
ordinary differential equations
\begin{equation}
\label{eq:l96}
\deriv{x_i}t=(x_{i-1}+\alpha\beta)(x_{i+1}-x_{i-2})-\beta
x_i+\beta^2(F-\alpha),
\end{equation}
with periodic boundary conditions $x_0=x_N$, with the total number of
variables $N=20$. The parameter $F>0$ provides constant forcing, and
the constant scaling parameters $\alpha>0$ and $\beta>0$ are chosen so
that, for given $F$, the statistical mean state $\langle
x_i\rangle=0$, and the statistical variance $\langle x_i^2\rangle=1$
(for details, see \cite{MajAbrGro}). For computation of the
correlation functions $\BS c_1(\tau)$ and $\BS c_2(\tau,s)$ from
\eqref{eq:corr}, we integrate the Lorenz 96 model in \eqref{eq:l96}
using the standard 4th order Runge-Kutta method, with the time
discretization step $\Delta t=0.01$ and the finite time averaging
window $T=10^6$ time units. The details of numerical discretization of
\eqref{eq:corr} are given in Appendix \ref{sec:app2}. For the
corresponding perturbed system, a small constant forcing $p$ is added
at a single node $x_i$; due to the statistical translational
invariance of \eqref{eq:l96}, the number $i$ of the node is
irrelevant. For the same reason, the entries of $\BS r(t)$ are
identical, $r_i(t)=r(t)$. Thus, the linear response approximation of
the Lyapunov exponent here is given by
\begin{equation}
\Delta\lambda\approx r(t_0)p.
\end{equation}
For the test, we pick three different values of $F=5,6$ and $8$, which
correspond to the dynamical regimes of low, moderate and strong chaos
and mixing, respectively. The standard time-lagged autocorrelation
functions $\langle x_i(t)x_i(t+\tau)\rangle$ (which are also identical
across different $i$ due to translational invariance of
\eqref{eq:l96}) for these regimes are shown in Figure
\ref{fig:l96_corr} as functions of the time lag. Observe that the
initial time scales of decorrelation are identical (the correlation
functions are almost the same for short lags), which is the effect of
rescaling by parameters $\alpha$ and $\beta$ in \eqref{eq:l96}. For
longer correlation lags we, however, can observe better mixing (more
rapid decay of lag-correlations) for the regimes $F=6,8$. The
corresponding values of the Lyapunov exponent for these regimes are
$\lambda=0.2265,0.3024$ and $0.4253$, respectively. In Figures
\ref{fig:l96_F5}--\ref{fig:l96_F8} we demonstrate the computed
finite-time response operators $r(t)$ for the regimes of the Lorenz 96
model with $F=5,6,8$, as well as the actual directly perturbed and
measured responses of the Lyapunov exponents for the perturbed system
and corresponding finite-time linear response approximations. We
observe that the temporal behavior of $r(t)$ generally consists of
three stages:
\begin{figure}
\includegraphics[width=\textwidth]{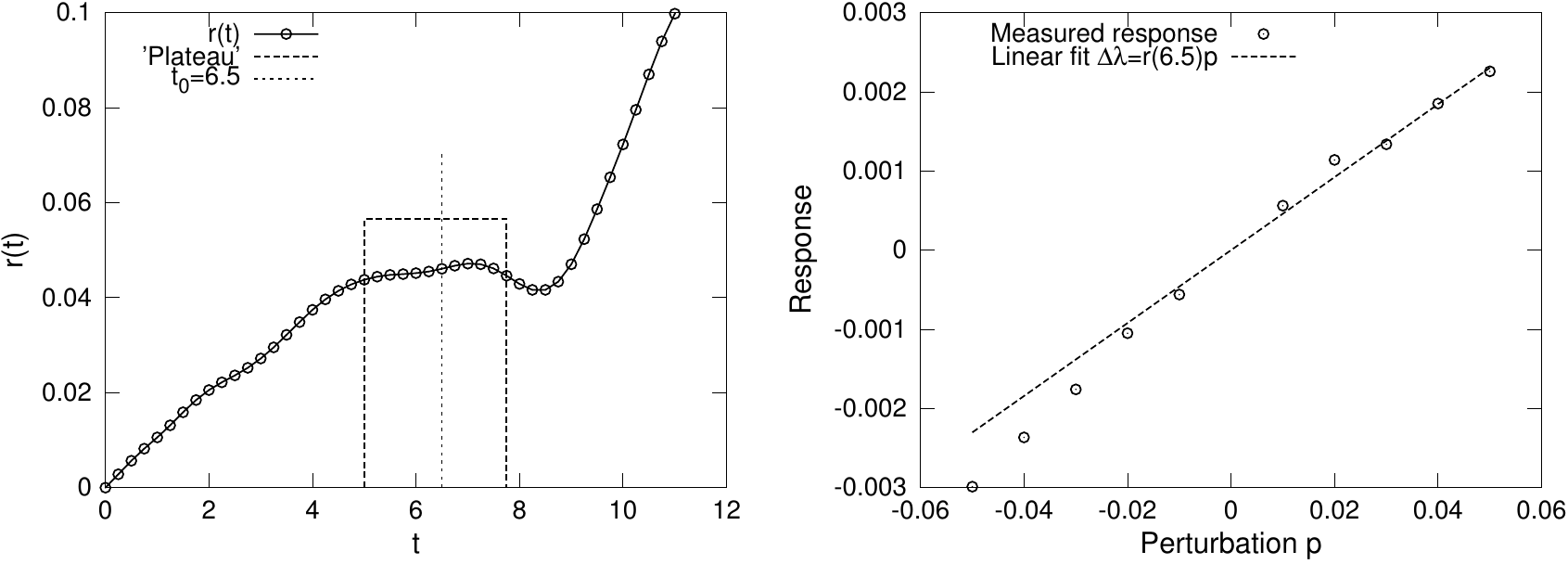}
\caption{Linear response of the Lyapunov exponent for the 20-variable
  Lorenz 96 model \eqref{eq:l96} with $F=6$ ($\lambda=0.3024$). Left:
  the plot $\BS r(t)$ as a function of finite response time, and the
  ``plateau'' of saturated response before manifestation of numerical
  instability. The time $t_0=6.5$ is chosen arbitrarily from within
  the ``plateau'' window. Right: the measured response of the Lyapunov
  exponent (via direct perturbations) vs the linear fit $\BS
  r(t_0)\cdot\BS p$.}
\label{fig:l96_F6}
\end{figure}
\begin{enumerate}
\item The initial growth stage, since $r(t)$ always starts at zero
  response for $t=0$.
\item The ``plateau'' stage, where the response has grown close to its
  equilibrium value. This stage should be the best approximation to
  the actual response of the Lyapunov exponent.
\item The blow-up stage, where the numerical instability in $\BS
  T_{\BS x}^t$ manifests itself. This stage is characterized by
  irregular oscillations and further growth of $r(t)$.
\end{enumerate}
The ``plateau'' stage of the linear response $r(t)$ is identified on
each plot, and the times $t_0$ for the finite-time response
approximations are chosen from within the plateau stages for
corresponding plots. It is not known at present whether the Lorenz 96
model is structurally stable in the displayed regimes (most likely
not), yet, structural instability does not appear to manifest itself
in the directly perturbed response significantly (that is, the circles
in the right-hand plots in Figures \ref{fig:l96_F5}--\ref{fig:l96_F8}
appear to lie on a ``smooth'' curve, rather than on some kind of a
discontinuous set).  Observe that the actual response of the Lyapunov
exponent to the perturbation is generally nonlinear, however, there is
a range of linearity around the unperturbed state, which appears to
depend on the dynamical regime. It is the most narrow in the weakly
chaotic regime ($F=5$), where the smallest perturbations used, $p=\pm
0.01$ yield linear response, however, larger perturbations ($p=\pm
0.02$ and above) cause a distinctly nonlinear, ``parabolic'' shape of
the response plot. As the dynamical regime becomes more chaotic and
mixing, the range of linearity becomes extended ($p=\pm 0.02$ for
$F=6$, and $p=\pm 0.03$ for $F=8$). In the range of linearity, the
finite-time linear response approximation yields a good fit to the
perturbed response of the Lyapunov exponent for each dynamical regime,
although an ``undershot'' (an insufficiently steep slope) can be seen
for the moderately chaotic and mixing dynamical regime $F=6$. Thus, it
appears that the weak chaos in the dynamics sets the natural bound of
applicability of the linearized formula, due to inherent nonlinearity
of the perturbed response of the Lyapunov exponent. The additional
observed effect is that the ``plateau'' of the response approximation
formula \eqref{eq:corr} becomes more narrow with the increase in
chaos, due to the fact that the numerical instability in the
exponentially growing tangent map manifests itself earlier in time.
This naturally leads to the speculation that, for strongly turbulent
dynamical regimes, the numerical instability could occur even before
the initial response growth stage is completed, thus leading to the
absence of a discernible ``plateau'', thus setting another bound of
practical applicability of the linear response formula in
\eqref{eq:corr}. Overall, for the regimes considered, the key
``plateau'' stages in each response operator $r(t)$ are clearly
identifiable by sight in each plot, and the linear response fits,
provided by the finite time linear approximation from within these
``plateau'' stages, seem to be adequate approximations to the directly
perturbed responses of the Lyapunov exponent.

\section{Conclusions}
\label{sec:conclusions}

In the current work we develop a linear approximation to the response
of the Lyapunov exponent of a nonlinear chaotic dynamical system to a
small constant perturbation. This approximation is computed from a
long-term trajectory of the corresponding unperturbed system. The
approximation is based on the fluctuation-dissipation theorem approach
\cite{Abr5,Abr6,Abr7,AbrMaj4,AbrMaj5,AbrMaj6,Rue2} using the tangent
map of the underlying chaotic dynamical system. We numerically test
the new approach using the rescaled Lorenz 96 model
\cite{Abr5,Abr6,Abr7,AbrMaj4,AbrMaj5,AbrMaj6,Abr8,Abr9,Abr10,Abr11,MajAbrGro,Lor,LorEma}
with 20 variables, in three dynamical regimes: weakly, moderately, and
strongly chaotic.  We show that, despite the fact that an inherent
numerical instability due to exponentially growing in time tangent map
renders the formal infinite time limit infeasible for practical
computation, a finite-time linear response formula adequately
approximates the actual perturbed values of the Lyapunov exponent in
their range of linearity, for the regimes considered. We also observe
that, for the same range of perturbations, the nonlinearity of the
perturbed response is strongest in the weakly chaotic regime, and,
vice-versa, weakest in the strongly chaotic regime. Also, the
``plateau'' stage of the linear response approximation between the
initial growth and numerical instability, which is crucial for
choosing the correct finite response time, appears to shrink when the
dynamical regime is strongly chaotic. It leads to the speculation that
for strongly turbulent regimes the developed method could be rendered
inapplicable due to the complete absence of the ``plateau''
stage. Thus, the range of practical applicability of the method
appears to be limited by the response nonlinearity on the weak chaos
side, and by the rapidly developing numerical instability in the
tangent map on the strong chaos side.

In the future work, we plan to investigate the response of the
Lyapunov exponent to linear vector field perturbations. Under special
interest are conservative systems whose solutions preserve a quadratic
energy (such as the truncated Burgers-Hopf system, the Kruskal-Zabusky
system, and unforced, undamped Lorenz 96 system
\cite{AbrKovMaj,AbrMaj,AbrMaj2,AbrMaj3}). In such systems, a
skew-symmetric (in the energy metric) linear perturbation vector field
will preserve the solutions on the same constant energy surface, at
the same time affecting chaos and turbulence of the system. An
interesting problem would be the maximization of the Lyapunov exponent
under linear skew-symmetric perturbation for given nonlinear
conservative dynamics, constrained to a fixed constant energy surface.

{\bf Acknowledgments.} This work was supported by the National Science
Foundation CAREER grant DMS-0845760.


\appendix
\section{Details on derivation}
\label{sec:app1}

The chain rule, applied to $L_{\BS p}$ in \eqref{eq:lin_rel}, yields
\begin{equation}
\begin{split}
\left.\parderiv{L_{\BS p}}{\BS p}\right|_{\BS p=\BS
  0}=\left.\parderiv{L_{\BS p}}{\phi_{\BS p}^t\phi^{-t}\BS
  x}\right|_{\phi_{\BS p}^t\phi^{-t}\BS x=\BS
  x}\left.\parderiv{\phi_{\BS p}^t\phi^{-t}\BS x}{\BS p}\right|_{\BS
  p=\BS 0}+\\+\left.\parderiv{L_{\BS p}}{\psi_{\BS p,\phi^{-t}\BS
    x}^t\psi_{\BS x}^{-t}\BS w}\right|_{\psi_{\BS p,\phi^{-t}\BS
    x}^t\psi_{\BS x}^{-t}\BS w=\BS
  w}\left.\parderiv{\psi_{\BS p,\phi^{-t}\BS x}^t\psi_{\BS x}^{-t}\BS
  w}{\BS p}\right|_{\BS p=\BS 0},
\end{split}
\end{equation}
where the corresponding partial derivatives are given by
\begin{subequations}
\begin{equation}
\left.\parderiv{L_{\BS p}}{\phi_{\BS p}^t\phi^{-t}\BS
  x}\right|_{\phi_{\BS p}^t\phi^{-t}\BS x=\BS x}=\BS w^TD^2\BS f(\BS
x)\BS w,
\end{equation}
\begin{equation}
\left.\parderiv{L_{\BS p}}{\psi_{\BS p,\phi^{-t}\BS x}^t\psi_{\BS
    x}^{-t}\BS w}\right|_{\psi_{\BS p,\phi^{-t}\BS x}^t\psi_{\BS
    x}^{-t}\BS w=\BS w}=\big(D\BS
  f(\BS x)+D\BS f(\BS x)^T\big)\BS w.
\end{equation}
\end{subequations}
Now we denote $\BS y=\phi^{-t}\BS x$ and compute $\partial\phi_{\BS
  p}^t\BS y/\partial\BS p|_{\BS p=\BS 0}$. Observe that $\phi_{\BS p}^t\BS
y$ satisfies
\begin{equation}
\parderiv{}t\phi_{\BS p}^t\BS y=\BS f(\phi_{\BS p}^t\BS y)+\BS p,
\end{equation}
where the differentiation on both sides with respect to $\BS p$ yields
the linear equation
\begin{equation}
\parderiv{}t\left(\parderiv{\phi_{\BS p}^t\BS y}{\BS p}\right)
=D\BS f(\phi_{\BS p}^t\BS y)\parderiv{\phi_{\BS p}^t\BS y}{\BS p}+\BS I,
\end{equation}
where $\BS I$ is the identity matrix. Upon replacing $\BS y$ back
with $\phi^{-t}\BS x$, the solution is given
\begin{equation}
\label{eq:dphi_dg}
\left.\parderiv{\phi_{\BS p}^t\phi^{-t}\BS x}{\BS p}\right|_{\BS p=0}=
\int_0^t\BS T_{\phi^{-\tau}\BS x}^\tau\,\dif\tau,
\end{equation}
where the tangent map $\BS T_{\BS p,\BS x}^t=\partial\phi_{\BS p}^t\BS
x/\partial\BS x$ is computed from the equation
\begin{equation}
\label{eq:tangent_map}
\parderiv{}t\BS T_{\BS p,\BS x}^t=D\BS f(\phi_{\BS p}^t\BS x)\BS
T_{\BS p,\BS x}^t, \qquad \BS T_{\BS p,\BS x}^0=\BS I.
\end{equation}
For $\partial(\psi_{\BS p,\phi^{-t}\BS x}^t\psi_{\BS x}^{-t}\BS
w)/\partial\BS p|_{\BS p=\BS 0}$, we again denote $\BS y=\phi^{-t}\BS
x$, and, additionally, $\BS z=\psi_{\BS x}^{-t}\BS w$, thus switching
to the computation of $\partial(\psi_{\BS p,\BS y}^t\BS z)/\partial\BS
g|_{\BS p=\BS 0}$. At this point, recall that $\BS w=\BS v/\|\BS v\|$,
where the tangent vector is given by $\BS v=\BS T_{\BS p,\BS y}^t\BS
z$. Using the chain rule of differentiation again, we obtain
\begin{equation}
\left.\parderiv{\psi_{\BS p,\BS y}^t\BS z}{\BS p}\right|_{\BS p=\BS
  0}=\parderiv{\BS w}{\BS v}\left.\parderiv{\BS T_{\BS p,\BS
    y}^t\BS z}{\BS p}\right|_{\BS p=\BS 0}.
\end{equation}
The term $\partial\BS T_{\BS p,\BS y}^t\BS z/\partial\BS p$ satisfies
\eqref{eq:tangent_map}, differentiated by $\BS p$ on both sides:
\begin{equation}
\begin{split}
\parderiv{}t\left(\parderiv{\BS T_{\BS p,\BS y}^t\BS z}{\BS p}\right)
=D\BS f(\phi_{\BS p}^t\BS y)\parderiv{\BS T_{\BS p,\BS y}^t\BS z}{\BS
  p}+\\+D^2\BS f(\phi_{\BS p}^t\BS y):\left(\BS T_{\BS p,\BS y}^t\BS
z\otimes\parderiv{\phi_{\BS p}^t\BS y}{\BS p}\right),
\end{split}
\end{equation}
where the Frobenius product ``$:$'' is computed over the two covariant
indices of $D^2\BS f$. The solution is given by
\begin{equation}
\label{eq:dT_dg_2}
\begin{split}
\left.\parderiv{\BS T_{\BS p,\BS x}^t\BS z}{\BS p}\right|_{\BS p=\BS
  0}=\int_0^t\BS T_{\phi^\tau\BS y}^{t-\tau} D^2\BS f(\phi^\tau\BS
y):\\:\bigg(\BS T_{\BS y}^\tau\BS z\otimes\int_0^\tau\BS T_{\phi^s\BS
  y}^{\tau-s}\dif s\bigg)\dif\tau,
\end{split}
\end{equation}
where we took into account \eqref{eq:dphi_dg}. Now, observe that
\begin{equation}
\BS T_{\BS y}^\tau\BS z = \|\BS v\|\BS T_{\phi^{-t}\BS x}^\tau\BS T_{\BS
  x}^{-t}\BS w=\|\BS v\|\BS T_{\BS x}^{\tau-t}\BS w,
\end{equation}
which, upon substitution into \eqref{eq:dT_dg_2} and rearrangement of
dummy variables of integration, yields
\begin{equation}
\label{eq:dT_dg}
\begin{split}
\left.\parderiv{\BS T_{\BS p,\BS x}^t\BS z}{\BS p}\right|_{\BS p=\BS
  0}=\|\BS v\|\int_0^t\BS T_{\phi^{-\tau}\BS x}^\tau D^2\BS f(\phi^{-\tau}\BS
x):\\:\bigg(\BS T_{\BS x}^{-\tau}\BS w\otimes\int_{\tau}^t\BS
T_{\phi^{-s}\BS x}^{s-\tau}\dif s\bigg)\dif\tau.
\end{split}
\end{equation}
Next we compute
\begin{equation}
\parderiv{\BS w}{\BS v}=\frac 1{\|\BS v\|}\big(\BS I-\BS w\BS
w^T\big),
\end{equation}
which, together with \eqref{eq:dT_dg}, yields
\begin{equation}
\begin{split}
\left.\parderiv{\psi_{\BS p,\phi^{-t}\BS x}^t\psi_{\BS x}^{-t}\BS
  w}{\BS p}\right|_{\BS p=0}= \big(\BS I-\BS w\BS
w^T\big)\int_0^t\BS T_{\phi^{-\tau}\BS x}^\tau \times\\\times D^2\BS
f(\phi^{-\tau}\BS x):\bigg(\BS T_{\BS x}^{-\tau}\BS
w\otimes\int_{\tau}^t\BS T_{\phi^{-s}\BS x}^{s-\tau}\dif
s\bigg)\dif\tau.
\end{split}
\end{equation}
Combining the computed terms together under \eqref{eq:lin_rel}, we
obtain
\begin{subequations}
\label{eq:average}
\begin{equation}
\lambda^*-\lambda\approx\lim_{t\to\infty}\BS r(t)\cdot\BS p,
\end{equation}
\begin{equation}
\BS r(t)=\int_0^t\BS c_1(\tau)\dif\tau+
\int_0^t\dif\tau\int_\tau^t\BS c_2(\tau,s)\dif s,
\end{equation}
\begin{equation}
\BS c_1(\tau)=\int\BS w^TD^2\BS f(\BS x):\big(\BS w\otimes\BS
T_{\phi^{-\tau}\BS x}^\tau\big)\dif\rho(\BS x,\BS w),
\end{equation}
\begin{equation}
\begin{split}
\BS c_2(\tau,s)=\int\BS w^T\left(D\BS f(\BS x)+D\BS f(\BS
x)^T\right)\times\\\times\big(\BS I-\BS w\BS w^T\big)\BS
T_{\phi^{-\tau}\BS x}^\tau D^2\BS f(\phi^{-\tau}\BS x):\\:\Big(\BS
T_{\BS x}^{-\tau}\BS w\otimes\BS T_{\phi^{-s}\BS
  x}^{s-\tau}\Big)\dif\rho(\BS x,\BS w).
\end{split}
\end{equation}
\end{subequations}
Invoking Birkhoff's ergodic theorem and replacing the measure averages
with long-term time averages, we obtain the correlation functions in
\eqref{eq:corr}.

\section{Details on discretization}
\label{sec:app2}

We discretize $\BS c_1(\tau)$ and $\BS c_2(\tau,s)$ from
\eqref{eq:corr} as follows: first, we assume that $\tau$ can assume a
finite range of values $\tau_m=\{0\ldots hm\ldots hM\}$, where $h$ is
the discretization step, and $M$ is a positive integer, with $hM$
bounding the range of possible values of $\tau$ from above. For each
value $hm$ of $\tau$, $s$ assumes the values $s_n=\{hm\ldots hn\ldots
hM\}$, so that the array of discretized values $(\tau_m,s_n)$ is
triangular. Then, the integrals over $\BS c_1(\tau_m)$ and $\BS
c_2(\tau_m,s_n)$ in \eqref{eq:corr} are computed using the standard
trapezoidal quadrature rule.

The discretized correlation functions $\BS c_1(\tau_m)$, $\BS
c_2(\tau_m,s_n)$ are computed as follows. Let us first discretize the
trajectory $\BS x(t)$, $\BS w(t)$ into the set of vectors $\BS x_k=\BS
x(t_k)$, $\BS w_k=\BS w(t_k)$, where $k$ is the discretization index
of time $t$. Let us also denote the incremental (forward by $h$)
tangent map $\BS T_{\BS x(t_k)}^h=\BS T_k$. Then, any tangent map $\BS
T_{\BS x(t_k)}^{hi}$ can be written as a product of the incremental
tangent maps \cite{Abr5,Abr6,Abr7}
\begin{equation}
\BS T_{\BS x(t_k)}^{hi}=\prod_{j=0}^{i-1}\BS T_{k+j}.
\end{equation}
Additionally, the backward (in time) tangent maps are readily
available as inverses of the corresponding forward tangent maps:
\begin{equation}
\BS T_k^{-h}=\big(\BS T_{k-1}^h\big)^{-1}.
\end{equation}
The incremental tangent maps $\BS T_k$ are obtained naturally by
solving the linear equation \eqref{eq:tangent_map} with $\BS p=0$
along the trajectory $\BS x(t)$ between $t_k$ and $t_k+h$. Then, the
time-averaged formulas for $\BS c_1(\tau)$ and $\BS c_2(\tau)$ in
\eqref{eq:corr} are expressed as the following discretized averages:
\begin{subequations}
\label{eq:discretized_corr}
\begin{equation}
\begin{split}
\BS c_1(\tau_m)=\lim_{K\to\infty}\frac 1K\sum_{k=M}^{K+M-1}\BS
w_k^TD^2\BS f(\BS x_k):\\:\bigg(\BS w_k\otimes\prod_{j=0}^{m-1}\BS
T_{k-m+j}\bigg),\quad 0\leq m\leq M,
\end{split}
\end{equation}
\begin{equation}
\begin{split}
\BS c_2(\tau_m,s_n)=\lim_{K\to\infty}\frac 1K\sum_{k=M}^{K+M-1}\BS
w_k^T\big(D\BS f(\BS x_k)^T+\\+D\BS f(\BS x_k)\big)\big(\BS
I-\BS w_k\BS w_k^T\big)\prod_{j=0}^{m-1}\BS T_{k-m+j}D^2\BS f(\BS
x_{k-m}):\\:\bigg(\prod_{j=0}^{m-1}\big(\BS T_{k-j}\big)^{-1}\BS
w_k\otimes\prod_{j=0}^{n-m}\BS T_{k-n+j}\bigg),\quad m\leq n\leq M.
\end{split}
\end{equation}
\end{subequations}
The formulas above are computed ``on-the-fly'' along with a
discretized long-term trajectory $\BS x_k$, $\BS w_k$, $\BS T_k$, with
the only memory requirement is that, as the discretized trajectory is
computed, the last $M$ incremental tangent maps $\BS T_k$ are
temporarily stored for computation, where $M$ is a fixed number.  Note
that $h$ above is not necessarily the discretization time step for the
4th order Runge-Kutta integrator of the unperturbed equation in
\eqref{eq:dyn_sys}; in particular, here we use $h=0.25$ time units
(with $M=60$, so that the response time spanned by $\BS r(t)$ is at
most $15$ time units), while the Runge-Kutta time step is $\Delta
t=0.01$ (that is, the averages in \eqref{eq:discretized_corr} are
updated once per $25$ Runge-Kutta steps). This allows to balance the
workload in a multithreaded implementation of the algorithm, where the
computation of $\BS x_k$, $\BS w_k$, $\BS T_k$, and the updates to
\eqref{eq:discretized_corr} are spread between different CPUs. The
time averaging window is $T=10^6$ time units, which sets $K=4\cdot
10^6$.

\end{document}